\documentclass[11pt]{amsart}
\usepackage{geometry}                
\newtheorem{theorem}{Theorem}
\newtheorem{corollary}[theorem]{Corollary}    
\newtheorem{definition}{Definition}
\newtheorem{proposition}{Proposition}
\newtheorem{lemma}{Lemma}

\usepackage{graphicx, color}
\usepackage{amssymb}
\usepackage{epstopdf}
\title{The Convolution Ring of Arithmetic Functions and Symmetric Polynomials}
\author{Huilan Li and Trueman MacHenry}

\begin{document}
\maketitle

\begin{abstract}
Inspired by Rearick (1968), we introduce two new operators, LOG and EXP. The LOG operates on generalized Fibonacci polynomials giving generalized Lucas polynomials. The EXP is the inverse of LOG. In particular, LOG takes a convolution product of generalized Fibonacci polynomials to a sum of generalized Lucas polynomials and EXP takes the sum to the convolution product. We use this structure to produce a theory of logarithms and exponentials within arithmetic functions giving another proof of the fact that the group of multiplicative functions under convolution product is isomorphic to the group of additive functions under addition. The hyperbolic trigonometric functions are constructed from the EXP operator, again, in the usual way. The usual hyperbolic trigonometric identities hold.  We exhibit new structure and identities in the isobaric ring. Given a monic polynomial, its infinite companion matrix can be embedded in the group of weighted isobaric polynomials. The derivative of the monic polynomial and its companion matrix give us the different matrix and the infinite different matrix. The determinant of the different matrix is the discriminate of the monic polynomial up to sign. In fact, the LOG operating on the infinite companion matrix is the infinite different matrix. We prove that an arithmetic function is locally representable if an only if it is a multiplicative function. An arithmetic function is both locally and globally representable if it is trivially globally represented.  
\end{abstract}





\section{Introduction} 
\label{sec:Intro} 

The relation between the ring of symmetric polynomials and the unique factorization domain of arithmetic functions with the Dirichlet  (convolution) product has been introduced and investigated starting  with MacHenry \cite{TM} and continuing with MacHenry and Tudose \cite{TM2,MT},  and MacHenry and Wong \cite{MW,MW2}.  In these papers it is shown that the convolution group of multiplicative arithmetic functions is represented locally by evalutions of Schur-hook polynomials in the ring of symmetric polynomials.   For this purpose,  it turns out to be convenient to write the homogeneous symmetric polynomials on the Elementary Symmetric Polynomial basis. When this is done the ring of symmetric polynomials is called the \emph{isobaric ring} and the homogeneous polynomials, \emph{isobaric polynomials} (Section \ref{sec:isobaric}).  

We recall one of the fundamental theorems of symmetric polynomial theory,  namely,  given a monic polynomial,  its coefficients, excluding the leading one,  are elementary symmetric polynomials of its roots up to sign. If the polynomial is of degree   $k$, then we denote by  $t_j,\ j=1,\dots,k$, the negative of the coefficient of the monomial of degree $k-j$, which is an elementary symmetric polynomial up to sign.  An isomorphism of the isobaric ring to the ring of symmetric polynomials is implicitly defined by $t_{j} \leftrightarrow (-1)^{j+1} e_j$.  
The monic polynomial in question then can be written as $$X^k-t_1X^{k-1}-\cdots-t_jX^{k-j}-\cdots-t_k,$$and is called a \emph{core polynomial} (Section \ref{sec:Core}). This polynomial will play an important role in what follows.  

We recall that \emph{arithmetic functions} are maps from the set of natural numbers together with $0$ to the complex numbers $\mathbb{C}$.  We denote the domain of such functions by $\mathbb{N}$. Frequently,  the map will be restricted to a subring of $\mathbb{C}$.  The arithmetic functions are a commutative ring with identity under addition and ordinary multiplication.  If instead of ordinary multiplication,  the Dirichlet (convolution) product is used as the product,  then with ordinary addition,  the structure is still a commutative ring with identity, denoted by $\mathbb{\mathcal{A}}$ \cite{HS}. In fact,  $\mathbb{\mathcal{A}}$ is a unique factorization domain, Cashwell and Everett  \cite{CE}.
The \emph{multiplicative} arithmetic functions $\mathcal{M}$ are those functions $f$ in $\mathcal{A}$
such that $f(mn)= f(m)f(n)$ whenever $(m,n) = 1$ \cite{PJM}.  This is equivalent to saying that they are functions which are determined by their values on prime powers at each prime $p$.  We say that they are determined \emph{locally}.  In a similar way, a  function $f$ in $\mathcal{A}$ is \emph{additive} if whenever $(m,n) = 1$ then  $f(mn) = f(m) + f(n)$, and denote the set of such functions by  $\mathcal{S }$.  The multiplicative functions $\mathcal{M}$ under convolution belong to the group of units of the unique factorization domain $\mathcal{A}$.  It is also the case that the elements in  $\mathcal{S}$ form a group,  an additive group.  Moreover,  $\mathcal{S}$ is isomorphic to $\mathcal{M}$,  Rearick  \cite{R}.  It is also the case that any element in $\mathcal{M}$ is locally linearly recursive, that is, for each prime $p$,  the values of the function at the powers of $p$ can be computed linearly recursively,  where the recursion parameters are either finite or infinite in cardinality (Section \ref{sec:MAF}). 

In 1968, Rearick \cite{R,R2} introduced the notions of ``logarithm" and ``exponential" functions of arithmetic functions, and used them to prove the isomorphism of $\mathcal{M}$ and $\mathcal{S}$, along with a number of other isomorphism of interesting subgroups in $\mathcal{A}$.  This result was reproved by Jesse Elliott in 2008 \cite{JE} along with a number of other similar results. In 1968 \cite{R,R2}, Rearick also introduced what he called the trigonometry of numbers by using the properties of the exponential function in the usual way to define ``Sine'' and ``Cosine".  In the recent paper of Elliott \cite{JE},  he does not refer to the inverse operator of the logarithm as an exponential function, but uses his rather more economical, but equivalent definition of logarithm to define its inverse operator. Rearick did not have an explicit expression for his exponential operator,  but showed its existence and properties inductively. 

In  \cite{MW} a linearly recursive sequence of isobaric (that is, homogeneous symmetric) polynomials, the Generalized Fibonacci Polynomials (GFP) \cite{TM2},  was exhibited whose images under evaluation (in some suitable ring) yield exactly the set of multiplicative functions $\mathcal{M}$,  and another linearly recursive sequence of isobaric polynomials, the Generalized Lucas Polynomials (GLP), whose images yield the additive arithmetic functions $\mathcal{S}$.  These same sequences of polynomials were also employed to give useful representations of certain algebraic number fields \cite{MW2}.   In Section \ref{sec:LE},  following the work of Rearick \cite{R,R2}, we define two operators,  $\mathcal{L}$ and $\mathcal{E}$, which we have also called LOG and EXP, and which are mutually inverse. We are able to represent each of these operators explicitly as matrices, thus as linear operators.  They have the very pleasant property that $\mathcal{L}$ is an isobaric degree-preserving map of the sequence GFP onto the sequence GLP, and hence, $\mathcal{E}$ maps the sequence GLP onto the sequence GFP, again preserving isobaric degree.  In \cite{MW}, it was shown that $\mathcal{M}$ is represented by suitable evaluations of the polynomials in the sequence GFP. It was  known in \cite{MW} that GFP and GLP are related by an isobaric degree-preserving bijection implying a proof of the isomorphism of the two groups $\mathcal{M}$ and $\mathcal{S}$. In this paper, the operator $\mathcal{L}$ exhibits this isomorphism and identifies $\mathcal{S}$ as the additive functions of $\mathcal{A}$. The operators $\mathcal{L}$ and $\mathcal{E}$ also  give new information about the structure of the ring of symmetric polynomials.

The $k \times k$ \emph{companion matrix} $A$ (Section \ref{sec:Core}) of a core polynomial $\mathcal{C}$ of degree $k$ is  one of the chief tools in this paper. $A$ is invertible if and only if the constant term of  $\mathcal{C}$ is non-zero.  $A$ can be extended  to a $\infty \times k$-matrix, the \emph{infinite companion matrix} $A^\infty$, by appending the orbit of its row vectors under the action of $A$. If $A$ is invertible, it can be extended both from the top and from the bottom.  When  $\mathcal{C}$ is a \emph{generic core}  polynomial,  that is,  when its coefficients are indeterminates, then the entries of $A^\infty$ appear as a subset $\mathcal{W}$ in the isobaric ring.  $\mathcal{W}$ is in fact an abelian group under the convolution product.  Its entries are signed Schur-hook polynomials,  the right hand column is the GFP sequence, and the traces of each $k \times k$-block yield the GLP sequence.  
 
In  Section \ref{sec:Core} the \emph{different matrix}  $D$ is discussed.  This matrix is generated  by constructing the $A$-orbit of a vector $d_k$ derived by differentiating  $\mathcal{C}$. It turns out that the vector $ d_k$ is the same as the vector $l_n$ when $ n= k$, where $l_n$ is the vector induced by $\mathcal{L}$.  From this we deduce that the LOG of $A^\infty$ is element-wise $D^\infty$ (Section \ref{sec:F-LE}), where $D^\infty$ is the complete $A$-orbit of $d_k$.  The different matrix $D$ appears as a $k \times k$-matrix  inside the infinite different  matrix $D^\infty$,  and the determinant of $D$ gives the discriminant with respect to the core polynomial.  The right hand column of $D^\infty $ is the GLP sequence.

Again, following Rearick \cite{R,R2},  we have used the operator $\mathcal{E}$ to define ``hyperbolic trigonometric functions" denoted by $S$ and $C$, and have shown that they give the usual hyperbolic trigonometric identities (Section \ref{sec:Trig}), providing some new identities in the ring of symmetric functions. In Subsection \ref{sec:App} we calculate these trigonometry functions for the LOG of some of the arithmetic functions that have served as examples.

This paper concentrates on the consequences of the representation of the convolution ring of arithmetic functions in the isobaric ring.  This point-of-view also throws light on the structure of the isobaric ring itself, and, in particular, gives rise to a ``trignometry'' and a ``geometry'' of isobaric polynomials. We shall discuss such generalizations in a forthcoming paper. 

In Section \ref{sec:Rep} we distingquish between arithmetic functions that are \emph{locally representable} and arithmetic functions that are \emph{globally representable}. This is essentially a distinction between multiplicative arithmetic functions,  which can always be represented locally,  and arithmetic functions, which can be represented by GFP but are not multiplicative.  We use a theorem of Rearick to show that this is essentially a dichotomy.  See REMARK in Section \ref{sec:Rep}.

In Section \ref{sec:F-LE}, we call two linear recursive sequences contained in the isobaric ring \emph{companions} if one is the image of the other under the operator $\mathcal{L}$ . The sequences GFP and GLP are companions, $G= \mathcal{L}(F)$ (Proposition \ref{mytwo}).  From the above remarks we see that the columns of $A^\infty$ and the columns of $D^\infty$ are companions. The results of this paper show that the operators $\mathcal{L}$ and $ \mathcal{E}$ are closely related to differentiation and integration, at least when operating on the GFP sequence and the GLP sequence.  For example, it is known that any partial derivative $\frac{\partial}{\partial {t_j}}$ of the GLP sequence gives the GFP sequence (Section \ref{sec:isobaric}).  It was pointed out above that the LOG takes $A^\infty$ to $D^\infty$. These facts have consequences for the representation of arithmetic functions \cite{MW} and for the representation of number fields \cite{MW2} in the isobaric ring;  they also give additional information about the structure of the isobaric ring, that is, about the ring of symmetric polynomials. 

Also we show that the Catalan sequence is a particularly interesting example of an arithmetic function which is globally representable but not locally representable;  but that it is a unit in the convolution ring (Subsection \ref{sec:CmpSeq}).  It is, in fact, its own convolution inverse and has as companion sequence,
\begin{equation*}
\Xi(n) =
{2n-1 \choose n}.
\end{equation*}

Finally, Section \ref{sec:Char} includes an application to the character theory of the symmetric group $S_n$, and to the P\'olya theory of counting.


\section{The Ring of Isobaric Polynomials} 
\label{sec:isobaric} 

\subsection{Isobaric polynomials} 
\label{iso-poly}

For a fixed $k$, an \emph{isobaric polynomial} of \emph{degree} $k$ and \emph{isobaric degree} $n$ is a polynomial of the form $$P_{k,n}(t_1,\ldots,t_k) = \sum_{\alpha \vdash n} C_{\alpha}  t_1^{\alpha_1}\cdots t_k^{\alpha_k},$$  where $\alpha = (\alpha_1,\ldots,\alpha_k)$ with $\alpha_j\in \mathbb{N}$ and $\sum_{j=1}^k j\alpha_j =n$ and $C_{\alpha}\in \mathbb{Z}$.

The condition $\sum_{j=1}^k j \alpha_j = n$  is equivalent to:  $(1^{\alpha_1},\ldots,k^{\alpha_k})$ is a partition of  $n$, whose  largest part is at most $k$, and we write this in the abbreviated,  and somewhat unorthodox form,  as $\mathbf{\alpha}   \vdash  n.$  (Note that, given $k$, the vector $\{\alpha_i\}$ is sufficient to reconstruct the partition.) In general, $C_{\alpha}$ can be chosen from any suitable ring. These isobaric polynomials form a graded commutative ring with identity under ordinary multiplication and addition of polynomials, graded by isobaric degree. This ring is naturally isomorphic to the ring of symmetric polynomials, where the isomorphism is given by the involution:
$$t_{j} \rightleftarrows (-1)^{j+1} e_j$$
with $e_j$ being the $j$-th elementary symmetric polynomial in $k$ variables (see \cite{ Macd, MT}).  This isomorphism associates the Complete Symmetric Polynomials $h_n$ in the ring of symmetric polynomials with the Generalized Fibonacci Polynomials (GFP) $F_{k,n}$ in the isobaric ring.  It also associates the Power Symmetric Polynomials $p_n$ with the Generalized Lucas Polynomials (GLP) $G_{k,n}$ in the isobaric ring (see \cite{MT}, and Macdonald \cite{Macd}).  The correspondence between complete symmetric polynomials and GFP and the correspondence between power symmetric polynomials and GLP can be shown inductively using the fact that GFP and GLP are linearly recursions of degree $k$; that is, 
$$F_{k,n} = t_1F_{k,n-1} +t_2F_{k,n-2} +\cdots+t_kF_{k,n-k}, \quad n-k\geqslant 0$$ and 
$$G_{k,n} = t_1G_{k,n-1} +t_2G_{k,n-2} +\cdots+t_kG_{k,n-k}, \quad n-k\geqslant 0$$
with $$F_{k,0}= 1$$
and $$G_{k,0}=k.$$
Note here $G_{k,0}=k$ varies as $n$ varies. For a deeper discussion of this variation see \cite{MT2} the discussion of the differential lattice and  the weighted isobaric polynomials.
It will be shown below that this initial condition arises in a natural way. 

The GFP and GLP can be written explicitly as follows:

 $$F_{k,n} = \sum_{\alpha \vdash n} 
  \left( \begin{array} {c}   | \alpha | \\  \alpha_1,\ldots,\alpha_k  \\   \end{array}
\right)  t_1^{\alpha_1}\cdots t_k^{\alpha_k},$$
 
 $$G_{k,n} =  \sum_{\alpha \vdash n} 
  \left( \begin{array} {c}  | \alpha | \\  \alpha_1,\ldots,\alpha_k  \\   \end{array}\right) \frac{n} {|\alpha|}  t_1^{\alpha_1}\cdots t_k^{\alpha_k},$$
  
\noindent  where $|\alpha| = \sum_{j=1}^k \alpha_j.$

\subsection{Weighted isobaric polynomials}
\label{WIP}

Note that when $k=2$ and $t_1=1, t_2=1,$ $F_{2,n}(1,1)$ and $G_{2,n}(1,1)$ are, respectively, the  sequence of Fibonacci numbers and the sequence of Lucas numbers. On the other hand,  
for a fixed $k$ both the GFP and the GLP are linearly recursive sequences indexed by   $n$,  and as we shall see below,  the indexing can, in general, be extended to the negative integers with  preservation of the linear recursion property (see  \cite{MT,MW}). All other  recursive sequences of isobaric polynomials are linear combinations of isobaric reflects of sequences of Schur-hook polynomials and form a (free) $\mathbb{Z}$-module of Weighted Isobaric Polynomials, denoted by $\mathcal{W}$, (see \cite{MT,MT2}).
 It is a remarkable fact that the only isobaric polynomials that can be elements in linearly recursive sequences of isobaric polynomials are those that occur in one of the sequences in $\mathcal{W}$ \cite[Theorem 3.4]{MT}. Every sequence in $\mathcal{W}$ can be presented in a closed form whose structure explains the term \textit{weighted}. Using the notation $ P_{\omega,k,n}$ to denote a polynomial in $\mathcal{W}$, where $\omega=  (\omega_1,\ldots,\omega_k)$ is the $weight$ vector, usually taken to be an integer vector, $k$ is the degree of the core polynomials and $n$ is the isobaric degree, we have
$$ P_{\omega,k,n} =  \sum_{\alpha \vdash n} 
\left(
\begin{array}{c}\textit{}
 |\alpha|   \\\alpha_1,\ldots,\alpha_k
\end{array}
\right)
\frac{\sum{\alpha_j\omega_j}}{|\alpha|}  t_1^{\alpha_1}\cdots t_k^{\alpha_k},$$
where $k$ and $\omega$ are fixed and $n$ varies. The linear recursion for an arbitrary weighted isobaric polynomial $P_{\omega,k,n}$ is given by
$$P_{\omega,k,n}=t_1P_{\omega,k,n-1}+t_2P_{\omega,k,n-2}+\cdots+t_kP_{\omega,k,n-k},\quad n-k\geqslant 0$$
with $P_{\omega,k,0}=\omega_k$.

Both the GFP and the GLP are weighted sequences,  the weight vectors being given by, respectively, $(1,1,\dots,1,\ldots)$ for the GFP,  and $(1,2,\ldots,j,\ldots)$ for the GLP. The columns of the infinite companion matrix (see Section \ref{sec:Core}) have weight vectors of the form \\$(0,\ldots,0,1,\ldots,1,\ldots)$, where the number of zeros is one less than the column number counting from the right. The Schur-hook polynomial indexed by $(n,1^r)$ has the weight vector $(0,\ldots,0,(-1)^r,\ldots,(-1)^r,\ldots)$ with $r$ zeros. The GLP are alternating sums of all of the Schur-hook polynomials of the same isobaric degree (see \cite{Macd,MT,MT2}), that is, the traces of the $k\times k$-blocks in $A^\infty$.

Moreover, it is important to note that a new variable $t_k$ will appear in $P_{\omega,k,n}$  for the first time when $n=k$. Thus, for fixed $\omega$ and $k$ the $P_{\omega,k,j}$ are the same for all $|j|< k$. We call this the \emph{conservation principle}.  We illustrate this with a listing of the first five GFP, taking the point of view that $k$ varies as $n$ varies:
\begin{itemize}
  \item $F_{k,0} = 1$
  \item $F_{k,1} = t_1$
  \item $F_{k,2} = t_1^2 + t_2$
  \item $F_{k,3} = t_1^3 + 2t_1t_2+t_3$
  \item $F_{k,4} = t_1^4 + 3t_1^2t_2 + t_2^2+2t_1t_3+t_4$
  \item $F_{k,5} = t_1^5 + 4t_1^3t_2 +  3t_1t_2^2 +3t_1^2t_3+2 t_2t_3+2t_1t_4+t_5$
\end{itemize}
and for the GLP we have
\begin{itemize}
   \item $G_{k,0} = k$
  \item $G_{k,1} = t_1$
  \item $G_{k,2} = t_1^2 + 2t_2$
  \item $G_{k,3} = t_1^3 + 3t_1t_2+3t_3$
  \item $G_{k,4} = t_1^4 + 4t_1^2t_2 + 2t_2^2+4t_1t_3+4t_4$
  \item $G_{k,5} = t_1^5 + 5t_1^3t_2 +  5t_1t_2^2 +5t_1^2t_3+5 t_2t_3+5t_1t_4+5t_5$
\end{itemize}
When the degree of the core polynomial is clear we shall write $F_n$ and $G_n$ for $F_{k,n}$ and $G_{k,n}$.

\begin{proposition}\label{myzero} 

$\displaystyle\frac{\partial G_{k,n}}{\partial t_j} =  n F_{k,n-j}.$

\end{proposition}

\proof \begin{eqnarray*}
\frac{\partial G_{k,n}}{\partial t_j}     & = &\frac{\partial }{\partial t_j}  \sum_{\alpha \vdash n} 
  \left( \begin{array} {c}  | \alpha | \\  \alpha_1,\ldots,\alpha_k  \\   \end{array}\right) \frac{n} {|\alpha|}  t_1^{\alpha_1}\cdots t_k^{\alpha_k}\\
 & = &  n \sum_{\alpha \vdash n,\ \alpha_j \neq0} \frac{|\alpha|!}{\alpha_1!\cdots\alpha_k!} \frac{\alpha_j}{|\alpha|} t_1^{\alpha_1}\cdots t_j^{\alpha_j-1}\cdots t_k^{\alpha_k}\\
 & = &  n \sum_{\alpha \vdash n,\ \alpha_j \neq0} \frac{(|\alpha|-1)!}{\alpha_1!\cdots(\alpha_j-1)!\cdots\alpha_k!} t_1^{\alpha_1}\cdots t_j^{\alpha_j-1}\cdots t_k^{\alpha_k}\\
  & = &  n \sum_{\beta \vdash n-j} \frac{|\beta|!}{\beta_1!\cdots\beta_k!} t_1^{\beta_1}\cdots t_k^{\beta_k}\\
 &=& n \sum_{\beta\vdash n-j} 
  \left( \begin{array} {c}  | \beta | \\  \beta_1,\ldots,\beta_k  \\   \end{array}\right) t_1^{\beta_1}\cdots t_k^{\beta_k}\\
  &=& n F_{k,n-j}. 
\end{eqnarray*}  \qed

\subsection{Convolution product}
\label{convolution}

Let $P$ and $ P' $ be in $\mathcal{W}$,  then the \emph{convolution product} $$P_n \ast P'_n =  \sum_{j=0} ^n P_jP'_{n-j}.$$

In particular, if $F$ and $ F' $ are GFP functions evaluated at $[t_1,\ldots,t_k]$ and $[t'_1,\ldots,t'_k]$, respectively,  then $$F_n \ast F'_n =  \sum_{j=0} ^n F_jF'_{n-j}.$$ 
This convolution product can be regarded as a product induced by the  co-product on complete symmetric polynomial $h_n \xrightarrow{\Delta}\sum_{j=1}^n h_j \otimes h_{n-j}$ in the Hopf algebra structure of the ring of symmetric polynomials in infinite variables. The convolution product can be applied to any pair of linear recursive sequences in $\mathcal{W}$.  The convolution product induces a group structure on evaluated GFP with $-{\bf t}_n$ being the convolution inverse of $F_n$, where $$-{\bf t}_n=\left\{\begin{array}{cc}1 & \hbox{ if } n=0, \\-t_n & \hbox{ if } n\geqslant1.\end{array}\right.$$


\section{The Companion Matrix}
\label{sec:Core} 

\subsection{The core polynomial and the infinite companion matrix}

\label{ICM}

The \textit{core polynomial}  $$X^k - t_1X^{k-1}-t_2X^{k-2}-\cdots-t_k$$  ``controls"  the Generalized Fibonacci and Lucas Sequences of degree $k$.  Precisely, we have the following consequences of specifying what we might call a \emph{generic core}, that is, a core with variables $t_j$ as coefficients. We shall abbreviate the core polynomial by the notation $[t_1,\ldots,t_k]$.  With each such core polynomial we have a \emph{companion matrix}
$$A = \left(\begin{array}{ccccc}0 & 1 &0& \cdots &0\\0&0&1&\cdots&0 \\\vdots&\vdots&\vdots&\ddots&\vdots\\0&0&0&\cdots&1\\t_k & t_{k-1}&t_{k-2} & \cdots&t_1\end{array}\right).$$
We let the companion matrix operate on its last row vector on the right,  and append the image vector to the companion matrix as a new  last row.  We repeat this process obtaining a matrix with infinitely many rows.  Note that  $A$ is invertible if and only if  $t_k \neq 0 $.  Assuming that   $t_k \neq 0 $ we can also extend the matrix from the top row upward by operating on the top row with $A^{-1}$, getting a doubly infinite matrix, that is, one with infinitely many rows in either direction  and $k$-columns.  We call this the \emph{infinite companion matrix}:

 {$A^\infty = \left(\begin{array}{cccc}\vdots&\ddots&\vdots&\vdots \\(-1)^{k-1}S_{(-2,1^{k-1})}&\cdots& -S_{(-2,1)} & S_{(-2)} \\(-1)^{k-1}S_{(-1,1^{k-1})} &\cdots&- S_{(-1,1)} & S_{(-1)} \\(-1)^{k-1}S_{(0,1^{k-1})}&\cdots& -S_{(0,1)} & S_{(0)} \\(-1)^{k-1}S_{(1,1^{k-1})}&\cdots & -S_{(1,1)} & S_{(1)} \\(-1)^{k-1}S_{(2,1^{k-1})} &\cdots& -S_{(2,1)} & S_{(2)} \\(-1)^{k-1}S_{(3,1^{k-1})}&\cdots & -S_{(3,1)} & S_{(3)} \\(-1)^{k-1}S_{(4,1^{k-1})}&\cdots & -S_{(4,1)} & S_{(4)}\\ \vdots&\ddots&\vdots&\vdots\end{array}\right)$ 
 $=((-1)^{k-j}S_{(n,1^{k-j })})_{\infty\times k}$.
 \vspace{0.5cm}

Loosely speaking, we can think of the core polynomial as a ``generating'' function for the GFP, the GLP, as well as all other linear recursive sequences of Schur-hook polynomials, as follows:

The right hand column of  $A^{\infty}$ is just the GFP sequence with terms of degree $k$,  which, in fact, has been extended to functions indexed negatively.  The $k \times k$ contiguous blocks in $A^{\infty}$ are just the powers, positive or negative, of the matrix $A$.  The set of traces of these powers  consists of the GLP of degree $k$.  Each of the $k$ columns is a linear recursive sequence with \textit{recursion parameters}  $t_1,\ldots,t_k$ whose elements are signed Schur-hook polynomials. For details see \cite{MW}.  

Let    $ V = \{v_i\}$  be any collection of $k$-dimensional vectors, where $v_{i+1} = v_i A$ then the following proposition is easily proved.

\begin{proposition}\cite{MW}\label{myone}  The set  $V$ is a linear recursion with recursion parameters  $t_1,\ldots,t_k$.
\end{proposition}\qed

\begin{theorem}\cite[Theorem 1]{MW} \label{mytheoremrecursive}  Every numerical recursive sequence can be obtained by evaluating a GFP sequence of degree $k$ for integer values of $t_1,\ldots,t_k$.
\end{theorem}\qed

\noindent The converse of this theorem is trivially true.

\subsection{The infinite different matrix} 
\label{IDM}

For a fixed $k$, let $\mathcal{C}'(X)$ be the derivative of $\mathcal{C}(X)$ then \begin{equation*}\mathcal{C}'(X) =kX^{k-1}-(k-1)t_1X^{k-2}-\cdots-t_{k-1}. \end{equation*} 

  \begin{definition}\label{mydefinition}  The coefficients of $\mathcal{C}'(X)$ determine a vector $$d_k=(-t_{k-1}, -2t_{k-2},\ldots,-(k-1)t_1,k),$$ which gives rise to a $k\times k$-matrix $D$ by constructing the orbit of $d_k$ under the operation of $A$ on the right. Thus $$D=\left(\begin{array}{c}d_kA^0 \\d_kA \\\vdots \\d_kA^{k-1}\end{array}\right),$$ called the \emph{different matrix}. 
\end{definition} 

$\det D=(-1)^{\frac{k(k-1)}{2}}\Delta$, where $\Delta$ is the discriminant of the core polynomial $\mathcal{C}(X)$. We call the complete orbit of the action of $A$ on $d_k$, the \emph{infinite different matrix} $D^\infty$.
 
\bigskip\noindent\textbf{REMARK.}
 The companion matrix generates linear recursions with the recursion parameters $(t_1,t_2,\ldots,t_k)$. Thus, the orbit of a vector  $v$ operated on by $A$, $vA^n$, gives columns that are linearly recursive.

\bigskip\noindent\textbf{EXAMPLE.} Let $k=3$,  then the different matrix is 
$$D = \left(\begin{array}{ccc}-t_2 & -2t_1 & 3 \\3t_3 & 2t_2 & t_1 \\t_1t_3  & 3t_3+t_1t_2 & t_1^2+2t_2
\end{array}\right).$$

Note that the right hand column of $D$ is a part of GLP. In fact, the right hand column of $D^\infty$ is the GLP sequence, see Proposition \ref{LAD}.


\section{Multiplicative Arithmetic Functions} 
\label{sec:MAF} 

\subsection{Multiplicative arithmetic functions}
\label{MAF}

An \emph{arithmetic function} $f$ is a map from $\mathbb{N}$ to $\mathbb{C}$.
It is \emph{multiplicative}  if $f(mn) = f (m) f(n)$ whenever $(m,n) = 1$.  It is \emph{completely multiplicative} if $f(mn) = f (m) f(n)$ for all $m$ and $n$,  and \emph{specially multiplicative} whenever it can be written as the convolution product of two completely multiplicative functions \cite{{ReS}}. 

 It is well known that the multiplicative functions are invertible with respect to the convolution product and that in particular they form a group under this product. In \cite{TM}  it is shown that the group generated by the completely multiplicative functions,  usually called the \emph{group} of \emph{rational multiplicative functions} \cite{{JR}}, is an uncountably generated free abelian group.  In this paper (and in \cite{MW}), completely multiplicative arithmetic functions are called \emph{degree 1} functions,  and specially multiplicative  arithmetic functions are called \emph{degree 2} functions.

Clearly, a function is multiplicative if and only if its values are determined by its values at prime powers. In the case of a degree 1 function (completely multiplicative),  it is determined by its values at primes.  In the case of a degree 2 function (specially multiplicative),  it is determined by its values at primes and primes squared.  We shall say  that multiplicative functions are defined \emph{locally}.  It may be that we need all of the prime powers to define a multiplicative function.  In any case the multiplicative functions are a subgroup of the group of units of the convolution ring $\mathcal{A}$. In \cite{MW} it is shown that multiplicative functions can be represented by evaluations of isobaric polynomials; in particular, by the polynomials in $\mathcal{W}$, or even more partiularly by the evaluations of the GFP.  We shall give examples of this in the next section,  and show the utility of this representation in the following sections.

The representation map takes a multiplcative function $f$ of degree $k$ (finite or infinite) to $F_{k,n}$, where the variables $t_j$ are evaluated at the convolution inverse of the function $f$.  It is known that the inverse of a rational arithmetic function is determined by a generating function of the form $X^k-t_1X^{k-1} -\cdots-t_k$.  This function serves as the core polynomial or, in the case of an infinite inverse, a power series,  (see \cite{MW} for a complete discussion of this problem),  which in turn determines a GFP.  It is the values of this linear recursive sequence of GFP which faithfully represents the $\mathcal{M}$.  We look at the application of this idea in the following section.

\subsection{Examples}
\label{Ex} 
The \emph{zeta function}  $\zeta$ and the \emph{M\"{o}bius function} $\mu$ are well known to be convolution inverses of one another in $\mathcal{M}$. $\zeta$ is a degree 1 function.  Suppose that what we know is that  $\zeta(n) = 1$ for all values of $n$,  so in particular, for powers of all primes $p$.  Let $F_n$ represent $\zeta$ at $p^n$.  Then $F_{k,0} =1 $ and $F_{k,1}= 1 = t_1$.  Since $F_2$ is assumed to have the value  1,  and $F_{k,2} = t_1^2+t_2 = 1$, so $t_2 =0$. It is then easily shown inductively, using the linear recursion,  that $t_j = 0$ whenever $j>1$.  From this we deduce that  the degree $k$ of $\zeta$ is 1.  The generating function for its inverse is the core polynomial $X-1$, we denote this by $core\ \zeta=[1]$. So the inverse of $\zeta$ is $-1$ at primes,  and $0$ at prime powers greater than 1, which is just the M\"obius function $\mu$. In general, given a multiplicative function $f$ we can use this technique to determine its degree and its inverse,  the core polynomial being the generating function for the inverse.  

We apply this argument to $\tau$ and $\sigma$, the functions that count, respectively, the number of divisors of $n$ and the sum of the divisors of $n$. We know that  $\tau(p^n) =n+1$ and $\sigma(p^n) = 1+p+\cdots+p^n$. Then for $\tau$, $t_1=2$, $t_2 = -1$, $t_j=0,\ j>2,$ that is, $core\ \tau=[2,-1]$, the degree is 2;  and for $\sigma$, $t_1 = p+1$, $t_2 = -p$, $t_j = 0,\ j>2$, that is, $core\ \sigma=[p+1,-p]$, the degree is also 2.  Conversely,  if we are given the core polynomial,  we know the inverse, the degree, and hence,  all values of the function.  

An interesting case is the Euler totient function: $\phi$.  The above argument gives us  $t_j = p-1 $ for all values of $j \geqslant1$.  In \cite{MW} we introduced the classification of multiplicative functions according to the cardinality of the non-zero part of the range and the range of the inverse.  There are four categories for the ordered pair
(range, inverse range): (infinite, finite), (finite, infinite), (infinite, infinite), (finite, finite). The $\zeta$ function is an example of the first category,  as are both of the multiplicative functions  $\tau$ and $ \sigma$,  while $\mu$ is an example of the second category, as are the convolution inverses of   $\tau$ and $ \sigma$.  The Euler totient function  $\phi$ is an example of the third category,  while the only multiplicative function in the last category is the convolution identity  \cite[Theorem 7(1)]{MW}.  

Consider the \emph{Catalan numbers}
 $\Gamma(n) = \frac{1}{n+1}{2n \choose n}$. The same procedure as in the previous examples identifies $t_j $ as $\Gamma(j-1)$, and in the representation by GFP above, $F_n = \Gamma(n)$.  However,  this representation is not local,  and the $\Gamma(n)$ is not multiplicative, although,  it is clearly a unit,  having a convolution inverse. This is an example which shows that the group of units of the convolution ring of arithmetic functions is larger than the convolution group of multiplicative functions.


\section{LOG and EXP Operators in the Isobaric Ring}
\label{sec:LE} 

\subsection{LOG and EXP}
\label{LandE} 
In two papers  by David Rearick \cite{R,R2},  the notions of ``Logarithm'' and ``Exponential'' operators of arithmetic functions were introduced.  These operators were inverses of one another. The Logarithm operator takes convolution products to sums in $\mathcal{A}$,  and the Exponential operator takes sums to convolution products.  He denoted these operators by $L$ and $E$. The Exponential operator was then used to define ``Sine'', ``Cosine'', and ``Tangent'' ($S,C,T$, respectively). Rearick used $L$  and  $E$  to show that certain groups of arithmetic functions were isomorphic.  To our knowledge,  there was no follow-up to these definitions and the isomorphism results in the later literature until 2008 \cite{JE}, where a different version of the operator $L$ is introduced and the isomorphisms are reproved without reference to Rearick's work. We shall extend these definitions to the isobaric ring.

We define two operators $\mathcal{L}$ and $\mathcal{E}$, which we also refer to as LOG and EXP. $\mathcal{L}$ and $\mathcal{E}$ are graded operators, that is, respect isobaric degree.

\begin{definition} \label{mydefinitionP} For a fixed $k$ and $n\geqslant 1$, 
$$\mathcal{L}(P_{n}) = -t_{n-1} P_{1} - 2t_{n-2}P_{2} -\cdots-(n-1)t_1P_{n-1} + nP_{n},$$
where $P_{n}$ is in $\mathcal{W}$ and $t_i=0$ for $i>k$.
\end{definition} 

\begin{proposition}For a fixed $k$ and $n\geqslant 1$, $$\mathcal{L}(P_{n}) =  - t_{k-1}P_{n-k+1} -\cdots-(k-1)t_1P_{n-1} + kP_{n}.$$
\end{proposition}\qed

In general for any $n$, we define
$$\mathcal{L}(P_{n}) =  - t_{k-1}P_{n-k+1} -\cdots-(k-1)t_1P_{n-1} + kP_{n}.$$
So $$\mathcal{L}(F_{k,0}) = k,$$  
where the core polynomial  for $F$ is given by  $[t_1,\ldots,t_k]$.

\begin{proposition}\label{mypropositionL:F---G} 
$$\mathcal{L}(F_n) = G_n,$$
\end{proposition}

\proof  Clearly $\mathcal{L}(F_{k,0}) = G_{k,0}$ and $\mathcal{L}(F_{k,1})= G_{k,1}$. 

Moreover, since $\mathcal{L}$ is a linear operator, by induction
\begin{eqnarray*}
\mathcal{L}(F_n) &=& \mathcal{L}(t_1F_{n-1}+t_2F_{n-2}+\cdots+t_nF_0)\\
&=& t_1\mathcal{L}(F_{n-1})+t_2\mathcal{L}(F_{n-2})+\cdots+t_n\mathcal{L}(F_0)\\
&=& t_1G_{n-1}+t_2G_{n-2}+\cdots+t_nG_0\\
&=& G_n.
\end{eqnarray*}  \qed

The operator $\mathcal{L}(P_n)$ can also be written as a matrix
$$L_n= \left(\begin{array}{ccccc}1 & 0 & 0 & \ldots & 0 \\-t_1 & 2 & 0 & \cdots & 0 \\-t_2 & -2t_1 & 3 & \cdots & 0 \\\cdots & \cdots & \cdots &\ddots & \vdots\\-t_{n-1} & -2t_{n-2} & -3t_{n-3} & \cdots & n\end{array}\right)$$
operating on the vector $(P_1,\ldots,P_n)$.

The matrix  
$$E_n=\left(\begin{array}{ccccc}1 & 0 & 0 & \cdots& 0 \\\frac{1}{2}F_1 & \frac{1}{2} & 0 & \cdots & 0 \\\frac{1}{3}F_2 & \frac{1}{3}F_1 & \frac{1}{3} & \cdots & 0 \\\vdots&  \vdots&  \vdots & \ddots& \vdots \\\frac{1}{n}F_ {n-1} & \frac{1}{n}F_{n-2} & \frac{1}{n}F_{n-3} & \cdots & \frac{1}{n}\end{array}\right) $$ is clearly the inverse of the matrix $L_n$,  thus  

\begin{definition} \label{mydefinitionE} For a fixed $k$,
$$\mathcal{E}(G_{k,0}) = 1$$  
$$\mathcal{E}(G_{k,n}) =\frac{1}{n} (F_{k,n-1} G_{k,1} +F_{k,n-2}G_{k,2} +\cdots+F_{k,1}G_{k,n-1} + G_{k,n}).$$
\end{definition}

\begin{proposition}\label{mytwo} 
$\mathcal{L}$ and  $\mathcal{E}$ are inverses of one another on $F$ and $G$, i.e.,
$$\mathcal{L}(F_n) = G_n,$$
$$\mathcal{E}(G_n) = F_n.$$
\end{proposition}\qed

From the discussion in Section \ref{sec:Core} and using Proposition \ref{mypropositionL:F---G} we see that the LOG of GFP is encoded in the infinite companion matrix $A^\infty$.

Define $\lambda F_n= nF_n$,  then a simple calculation shows that
 
 \begin{lemma}\label{mylemmalog} For $n\geqslant 1$,              
 \begin{equation*} \mathcal{L}(F_n) = \overline{F_n}  \ast  (\lambda F)_n.  \end{equation*}
 \end{lemma}
\qed

  \begin{lemma}\label{mylemmalogadd} For $n\geqslant 1$,           
  \begin{equation*}\lambda(P'_n \ast P''_n)=P'_n \ast \lambda P''_n + P''_n  \ast  \lambda P'_n. \end{equation*}
  \end{lemma}\qed
 
Let $F'$ and $F''$ be two GFP induced by $[t_1',\ldots,t_k']$ and $[t_1'',\ldots,t_k'']$, respectively. Then $F'\ast F''$ is also a GFP induce by some core polynomial of degree $n$. Meanwhile, $G'$ and $G''$ are induced by $[t_1',\ldots,t_k']$ and $[t_1'',\ldots,t_k'']$, respectively.

\begin{proposition}\label{myLOG} For $n\geqslant 1$,
$$\mathcal{L}(F_n' \ast F_n'') = \mathcal{L}(F_n') + \mathcal{L}(F_n''),$$
$$\mathcal{E}(G_n' + G_n'') = \mathcal{E}(G_n') \ast \mathcal{E}(G_n'').$$
\end{proposition}

\proof
By Lemma \ref{mylemmalog}, $\mathcal{L}(F_n' \ast F_n'') = \overline{F_n' \ast F_n''} \ast \lambda (F_n' \ast F_n'')$. Multiplying both sides of Lemma \ref{mylemmalogadd}  by $\overline{F_n' \ast F_n''}$ and using distributive property of inverse over $\ast$ gives
 \begin{eqnarray*}
  &&\overline{F_n' \ast F_n''} \ast \lambda (F_n' \ast F_n'') \\
  &=& \overline{F_n'} \ast \overline{F_n''} \ast  ({F_n'} \ast\lambda F_n''+{F_n''} \ast\lambda F_n')\\
  &=&\overline{F_n' } \ast \lambda F_n' + \overline{F_n''} \ast \lambda F_n''.
\end{eqnarray*}
  So that $$\mathcal{L}(F_n' \ast F_n'') = \mathcal{L}(F_n') + \mathcal{L}(F_n'').$$ 
 While $$\mathcal{E}(\mathcal{L}(F_n' \ast F_n''))=\mathcal{E}( \mathcal{L}(F_n') + \mathcal{L}(F_n'')).$$ 
 So \begin{equation*} F_n' \ast F_n''=\mathcal{E}( \mathcal{L}(F_n') + \mathcal{L}(F_n'')).\end{equation*}
 Thus \begin{equation*} \label{E_Formula}  \mathcal{E}(G_n') \ast \mathcal{E}(G_n'')  =   \mathcal{E}( G_n' + G_n'').  \end{equation*} \qed

By  $P^{\ast r} $ we mean the $r$-th convolution power of $P$, for $P\in\mathcal{W}$.

\begin{corollary}\label{myl=d} For $n\geqslant 1$,
 \begin{equation*}
\label{rLOGF}
\mathcal{L}(F_n^{\ast r}) = r\mathcal{L}(F_n).
\end{equation*}
\begin{equation*}
\label{ }
\mathcal{E}(rG_n)=\mathcal{E}(G_n)^{\ast r}.
\end{equation*}
 \end{corollary}
 \proof Let $r\in\mathbb{N}$, then using Proposition \ref{myLOG}  and induction we have
\begin{equation*}
\label{ } 
\mathcal{L}(F_n^{\ast r}) =  \sum\mathcal{L}(F_n)=r\mathcal{L}(F_n)\end{equation*}
and 
 \begin{equation*}
\label{ }
\mathcal{E}(rG_n) = \mathcal{E}(r\mathcal{L}(F_n))=\mathcal{E}(\mathcal{L}(F_n^{\ast r}))= F_n^{\ast r}=\mathcal{E}(G_n)^{\ast r}.
\end{equation*}\qed
\vspace{0.5cm}

In fact, in \cite[Theorem 5.7]{MT}, we embed the convolution group of $\mathcal{W}$ in its injective hull,  that is,  to each such polynomial we adjoin its $q$-th root for every rational number $q$, giving a divisible group, also see \cite{CG}.

\vspace{0.5cm}
\begin{proposition}\label{mypropositionL:Fisadditive} 
$\mathcal{L}$ is a local isomorphism of $\mathcal{M}$ to $\mathcal{S}$.  Thus $\mathcal{S}$ is represented locally by the GLP for $G_{k,n}$ having integer cores.
\end{proposition}

\proof This follows from Proposition \ref{mypropositionL:F---G} and Proposition \ref{mytwo}. \qed

\subsection{Examples}
\label{Ex1} 
Consider the multiplicative functions $\zeta (n) = 1$, for all $n$, $\tau(n) = \mathrm{card}\{d:d|n\}$ and $\sigma(n) = \sum_{d|n} d$.
$\zeta,\tau, \sigma  \in \mathcal{M}$.  \vspace{0.5cm}

\noindent(1)  $ core \,  \zeta  =  [1] $. Represent $\zeta$  locally at prime $p$ by $\zeta \xrightarrow{{\rm lr}_p}  \, _\zeta F$, where $ \,  _\zeta F_n = 1$.  

$$ \,_{\zeta}\mathcal{L}(\,  _\zeta F) = \, _\zeta G = \zeta = (1,1,\ldots).\;$$
In general, if $f\xrightarrow{{\rm lr}_p} \, _f F$ then $core\ f= core\ _fF$.
\vspace{0.5cm}

\noindent(2) $ core \; \tau = [2,-1]$. $\tau(p^n) = n+1$, in fact $\tau = \zeta \ast \zeta$. Represent $\tau$ locally at $p$ by   $\tau \xrightarrow{{\rm lr}_p} \, _ \tau  F,$ where $  _\tau F_n = n+1. $  $$\,_{\tau} \mathcal{L}(_\tau F)  =  \, _\tau G = (2,2,\ldots) = 2 \,_{\zeta}\mathcal{L} (_\zeta F)  = \,_{\zeta}\mathcal{L} (_\zeta F)+ \,_{\zeta} \mathcal{L}(_\zeta F).$$
\vspace{0.5cm}

\noindent(3) $ core \,\sigma = [p+1,-p]$.  $\sigma(p^n)=1+p+\cdots+p^n$, in fact $\sigma = \zeta \ast \zeta_1$, where $\zeta_1$ is the degree 1 function such that $\zeta_1 (m) = m$ for all $m$; in  particular, $\zeta_1(p^n) = p^n$ for all prime $p$ and all $n$.
 Represent $\sigma$ locally at $p$ by $\sigma \xrightarrow{{\rm lr}_p}    \, _\sigma F_n ,$ where $\, _\sigma F_n =  1+p+\cdots+p^n .$ 
$$ \,_{\sigma} \mathcal{L}(\, _\sigma F_n) = \, _\sigma G_n = p^n+1 =\, _\zeta\mathcal{L} ( \, _\zeta F_n)  + \,_{\zeta_1} \mathcal{L}( \, _{\zeta_1} F_n).$$
\vspace{0.5cm}

Note that $ core \,\sigma = [p+1,-p].$   If  we let $p$ be 1,  then we get the core for $\tau$.
\vspace{0.5cm}

All of the above arithmetic functions are in the  group of rational multiplicative functions.

\subsection{Some consequences}        
\label{conseq}

Above, we have an example of a degree 1 multiplicative function $\zeta$ on which  $\,_\zeta\mathcal{L}$ acts as the identity;  and hence so does  $\,_\zeta\mathcal{E}$.  This is,  in fact, a property of all degree 1 multiplicative functions.

\begin{proposition}\label{myidemLOG} 
If $f$ is a degree 1 multiplicative function,  then 
$$\,_f\mathcal{L}(\,_f F) = \,_fF,$$
$$\,_f\mathcal{E}(\,_f G) = \,_fG.$$
\end{proposition}

\proof Since degree 1 multiplicative functions have core polynomials $[t_1]$---that is,  all of the $t_j = 0,$ $ j>1$, it follows that the GFP that represents the arithmetic function is of the form $t_1^n$.  Hence $\,_f\mathcal{L}(\,_fF_n)=\,_fG_n=\,_fF_n$ and $\,_f\mathcal{E}(\,_fG_n)=\,_fF_n$, both statements in the proposition follow. \qed

Denote the convolution product of $\{f_j\}$ by $\prod^\ast_jf_j$. From this simple proposition and from Proposition \ref{myLOG} a very useful computational fact follows:

\begin{corollary}\label{mydegreeone} 
Let $f = \prod^\ast_j f_j$ with $f_j $ having degree 1 for all $j$, then $\,_f\mathcal{L}(\,_fF) = \sum_j \,_{f_j}\mathcal{L}(\,_{f_j}F)$. 
\end{corollary}

\proof Let $f$ be a rational multiplication function, then  $f$ is the convolution product of finitely many degree 1 multiplicative functions.  Thus,  the theorem follows by induction, using Proposition \ref{myLOG}. \qed

Our definition of the operators $\mathcal{L}$ and $\mathcal{E}$ is an adaptation of Rearick's 1968 definition \cite{R,R2}.

\begin{definition}\cite[Definition 1]{R} 
If $f$ is an arithmetic function, let
\begin{eqnarray*}
Lf(n)&=& \sum_{d|n}f(d)f^{-1}(n/d) \log d, \qquad\qquad\hbox{ if }n > 1,\\
Lf(1)&=& \log f(1).
\end{eqnarray*}
\end{definition}

We have adapted the Rearick logarithm \cite[Definition 1]{R}  for use with multiplicative functions by noting that  at a prime $p$, this definition becomes
  \begin{eqnarray*}
Lf(p^n)&=& \sum_{j=0}^nf(p^j)f^{-1}(p^{n-j}) \log p^j ,\qquad\qquad\hbox{ if }n\geqslant 1,\\
Lf(p^0)&=& \log f(p^0).
\end{eqnarray*}
Then, for each prime $p$, we let  $f \xrightarrow{{\rm lr}_p} F$, giving
\begin{eqnarray*}
\mathcal{L}(F_n) &=& \sum_{j=0}^nj F_j\overline{F_{n-j}}\log p, \qquad\qquad\hbox{ if }n \geqslant 1.
\end{eqnarray*}
Finally, we take $\log p$ base $p$ yielding
\begin{eqnarray*}
\mathcal{L}(F_n) &=& \sum_{j=0}^nj F_j\overline{F_{n-j}},  \qquad\qquad\hbox{ if }n \geqslant 1.
\end{eqnarray*}

\bigskip
\noindent\textbf{REMARK.} Rearick's logarithm differs from ours at each prime $p$ by $\log p$. In particular, Rearick's logarithm of $\zeta$ is the  Mangoldt function while our logarithm of a degree 1 function is just itself.

\section{Further Properties of $\mathcal{L}$ and $\mathcal{E}$}
\label{sec:F-LE} 

\subsection{Identities}        
\label{iden}

The following identities can be used to establish a useful alternate description of the LOG and EXP operators.

 \begin{lemma}\label{mylemmaFandG}  $G_1 = F_1$ \\
$G_2 = -F_1^2 + 2F_2 $\\
$G_3 = F_1^3 -3F_1F_2+3F_3$\\
$G_4 = -F_1^4+4F_1^2F_2-2F_2^2-4F_1F_3+4F_4$\\
and, in general,
\begin{equation}\label{GF}
G_n = \sum_{\alpha \vdash n} n(-1)^{|\alpha|+1} \left(\begin{array}{c}|\alpha| -1  \\\alpha_1,\ldots,\alpha_k\end{array}\right) F_1^{\alpha_1}\cdots F_k^{\alpha_k}.
\end{equation}\qed
\end{lemma}

Note that the general  expression (\ref{GF}) above is the definition of  $G$ with signs altered according to the rule $(-1)^{|\alpha|+1}$ and a parts lift to $F$,  where the term \emph{parts lift} is defined as follows:

Given the partition $(1^{\alpha_1}, 2^{\alpha_2},\ldots,k^{\alpha_k})$,  we call the expression  $$P^{\alpha}=P_1 ^{\alpha_1}P_2^{\alpha_2}\cdots P_k^{\alpha_k}$$ the \emph{parts lift} to $P\in\mathcal{W}$ of this partition.

In the inverse direction we have
\begin{corollary}\label{mycorollaryGtoF}

$F_1 = G_1$\\
$F_2 = \frac{1}{2!}G_1^2 + \frac{1}{2}G_2$\\
$F_3 = \frac{1}{3!}G_1^3 +\frac{1}{2}G_1G_2 + \frac{1}{3}G_3$\\
$F_4 = \frac{1}{4!}G_1^4 + \frac{1}{4}G_1^2G_2+\frac{1}{8}G_2^2+\frac{1}{3}G_1G_3+\frac{1}{4}G_4$  \\
and, in general, 
\begin{equation*}F_n = \sum_{\alpha \vdash n} \frac{1}{z(\alpha)}G_1^{\alpha_1}\cdots G_k^{\alpha_k}.\end{equation*}
\end{corollary}\qed

\begin{corollary}\label{mycorollaryGtoG} 
\begin{equation*}\mathcal{E}(G_n)=\sum_{\alpha \vdash n} \frac{1}{z(\alpha)}G_1^{\alpha_1}\cdots G_k^{\alpha_k}. \end{equation*}
\end{corollary} \qed

\subsection{The infinite companion matrix and the infinite different matrix} 
\label{ICMIDM}

For a fixed $k$, we shall mean by $\mathcal{L}(M)$ a matrix whose entries are $(\mathcal{L}m_{(i,j)})$, where$$\mathcal{L}m_{(i,j)}= - t_{k-1}m_{i-k+1,j} -\cdots-(k-1)t_1m_{i-1,j} + km_{i,j}.$$  

It is convenient at this point to introduce the following notation: We denote the vector induced by the operator $\mathcal{L}$ as $l_n$ where $l_n$  is the $n$-th row vector of the matrix $L_n$.  Thus,  the third row of the matrix $L_3$ yields the vector $l_3=(-t_2, -2t_1, 3)$.  

Observe that $\mathcal{L}$ preserves isobaric degree, for example, $\mathcal{L}(F_{k,n})=G_{k,n}$. On the other band, $d_k$ respects the grading of the isobaric ring induced by the degree of the core polynomial. In other words, the definition of $l_n$ is independent of $k$ and the definition of $d_k$ is independent of $n$. But we have the important remark:

\bigskip
\noindent\textbf{REMARK.}
$l _n= d_k$, when $n=k$.
 
 From this there follows a host of consequences: 

 \begin{proposition}\label{LAD}
 \begin{equation*}
\mathcal{L}(A^\infty)  = D^\infty.
\end{equation*} 
\end{proposition}\qed
 
In particular, \begin{equation*}
\mathcal{L}\{\hbox{right-hand column of }A^\infty \} = G,
\end{equation*}
that is
\begin{equation*}
\mathcal{L}(F_n) = G_n,
\end{equation*}
and
\begin{equation*}
\mathcal{L}((-1)^r{S_{(n,1^r)}} )= d_{(n+1,k-r)},
\end{equation*}
where $d_{(i,j)}$ is the $(i,j)$-th element of $D^\infty.$
\vspace{0.5cm}

Since we know that the operator $\mathcal{L}$ takes multiplicative arithmetic functions to additive arithmetic functions (see Section \ref{sec:LE}),  and since $\mathcal{W}$ represents multiplicative functions, we have the following:

 \begin{corollary}\label{myadd} 
 The columns of $D^\infty$ represent additive arithmetic functions.
 \end{corollary}

\noindent\textbf{REMARK.}  In \cite{MW2} it was shown that the rings, and in particular the number fields of the form $\mathbb{Q}[X]/ideal<\mathcal{C}(X)>$, are representable in $\mathcal{W}$. The results just stated in this paper sharpen the results of \cite{MW2}.  The ring of symmetric functions, and in particular $\mathcal{W}$ give a strong link between arithmetic number theory and algebraic number theory, and, incidentally, the theory of linear recursions.  A core polynomial immediately induces a set of multiplicative arithmetic functions, a set of additive arithmetic functions and a number ring, as well as a set of linear recursions,  all of which are represented inside the isobaric ring.  In this paper and in previous papers the coefficient ring has been restricted to the rationals and integers.  It is clear that this can be generalized considerably. This will be taken up in subsequent papers.

\subsection{Companion Sequences}    
\label{sec:CmpSeq} 

In \cite{L}, D. H. Lehmer introduced the idea of pairs of companion sequences. We will call two arithmetic functions \emph{companions} if they are closely related in the way that the Fibonacci sequence and the Lucas sequence are related. We make this notion of ``closely related'' precise in the following definition.
 
\begin{definition} \label{mydefinitioncompanion}
If an arithmetic function $f$ is representable then its logarithm $g$ is its \emph{companion}, or we say that the two arithmetic functions $f$ and $g$ are \emph{companions}.
\end{definition}

\noindent\textbf{REMARK.} If $f$ and $g$ are companions then they have the same core, and $f$ is represented in GFP and $g$ is represented in GLP.  Thus, $f$ is a multiplicative arithmetic function and $g$ an additive arithmetic function. Moreover, $f$ and $g$ are represented in the same evaluation of the infinite companion matrix $A^\infty$.
 
As a result of this discussion, we have already computed the companion sequences for $\tau, \sigma, \zeta$ and $\phi$.  Here we will exhibit another interesting examples. In Section \ref{Ex}, we mentioned the Catalan numbers
 $\Gamma(n)$.  This sequence is globally representable for $n\geqslant 1$, $\Gamma(n)\xrightarrow{{\rm gr}}   \,_\Gamma F_n $.  So its companion is given by $\,_\Gamma\mathcal{L}(\,_\Gamma F_n) = \,_\Gamma G_n$.  We shall exhibit this sequence precisely, 
but first we mention an interesting property of the representation of the Catalan sequence.  We remarked in Section \ref{Ex} that $\Gamma(n)$ is, what we call, an \emph{incestuous sequence}; namely, it is self-inverse, that is,  the sequence $t_j$ which constitutes the inverse of $\Gamma$ is again $ \Gamma$, 
\begin{equation*}
\label{ }
t_{j+1} = \Gamma(j).
\end{equation*}
An equivalent version of this fact is contained in the following

 \begin{proposition}\label{myglobnotloc}For $n\geqslant 1$,
 \begin{equation*}
\,_\Gamma F_n \ast \,_\Gamma F_n = \,_\Gamma F_{n+1}
\end{equation*}
\end{proposition}\qed

We use Proposition \ref{mypropositionL:F---G} to compute $\,_\Gamma\mathcal{L}(\,_\Gamma F_n)$.

 \begin{proposition}\label{myglobnotloc} For $n\geqslant 1$,
 \begin{equation*}
\,_\Gamma\mathcal{L}( \,_\Gamma F_{n}) =
{2n-1 \choose n}.
\end{equation*}
 \end{proposition}
 
\proof Since $\,_\Gamma\mathcal{L}(\,_\Gamma F_{n})= \,_\Gamma G_{n}$ and $\,_\Gamma G_{n}$ is a linear recursion, this proposition follows easily by induction. \qed
 
 Note that \begin{equation*}
\,_\Gamma G_{n} = \frac{n+1}{2} \Gamma(n)
\end{equation*}

Denote the sequence $\,_\Gamma G_{n}$ by $\Xi(n)$. $\Gamma(n)$ and $\Xi(n)$ are associated with Dyck paths, among other things  (see sequence A001700 in  \cite{online}).


\section{Character tables and P\'olya's counting theorem}
\label{sec:Char} 

As a result of Corollary \ref{mycorollaryGtoG},  the character theory of the finite symmetric group, as well as P\'olya's Counting Theorem, can be conveniently described in terms of the sequences GFP and GLP.  

First we recall that the entries in the infinite companion matrix in Section \ref{sec:Core} are Schur-hook polynomials,  and that these polynomials can be recovered as isobaric polynomials using the Jacoby-Trudi formula with entries from among the generalized Fibonacci polynomials.  This can be done in general for Schur polynomials$$S_\lambda= \det(F_{\lambda_i - i+j}) ,$$ where $\lambda$ is listed in weakly decreasing order.  For example the Schur polynomial $S_{(3,2)}$   induced by the partition $(3,2)$ in isobaric form  is given by 
\begin{eqnarray}
\label{Schur}
S_{(3,2)}& =& \det \left( \begin{array}{cc }
  F_3 & F_4    \\
  F_1  & F_2   
\end{array}\right) \nonumber\\
&=& F_{2}F_3 - F_1F_4 \nonumber\\
&=& t_1t_2^2-t_1^2t_3+t_2t_3-t_1t_4\nonumber
\end{eqnarray}

\noindent{\bf Character Algorithm.} \cite[p.114]{Macd}
 The irreducible character $\chi^{\lambda}$ of $S_n$ indexed by the partition $\lambda = (\lambda_1, ..., \lambda_r)$, written in weakly descreasing order, is given by
\begin{eqnarray}
ch(\chi^{\lambda}) &=&  \det(F_{\lambda_i - i+j}) \nonumber\\
\label{det}&= &\det (\mathcal{E}(G_{\lambda_i - i+j}))
\end{eqnarray}
Applying Corollary \ref{mycorollaryGtoG}, we get that the result of (\ref{det}) is a polynomial of the form 
\begin{equation*}
\label{ }
\sum_{\alpha \vdash n}c(\alpha)G^{\alpha}.
\end{equation*}
This expression can be mapped to the vector  $( \chi^{\lambda}(\alpha))_{\alpha\vdash n}$ by the mapping
$$ \sum_{\alpha \vdash n}c(\alpha)G^{\alpha} \rightarrow ( \chi^{\lambda}(\alpha))_{\alpha\vdash n}=(c(\alpha)z(\alpha))_{\alpha\vdash n}.$$

\vspace{0.5cm}

\noindent\textbf{EXAMPLE.}  Let $n=4$ and let the Schur polynomial be $S_{(3,1)}$,  then
\begin{eqnarray*}
S_{(3,1)}& =& \det \left( \begin{array}{cc}
   F_3   & F_4   \\
   F_0   & F_1  
\end{array} \right)\\
&=& F_1F_3-F_4\\
&=& \frac{1}{6} G_1^4 + \frac{1}{2} G_1^2G_2+ \frac{1}{3} G_1G_3 - \frac{1}{4!} G_1^4 - \frac{1}{4} G_1^2G_2 - \frac{1}{8} G_2^2 - \frac{1}{3} G_1G_3 - \frac{1}{4} G_4 \\
&= &\frac{1}{8} G_1^4 +\frac{1}{4} G_1^2G_2 -\frac{1}{8} G_2^2 -0 G_1G_3 - \frac{1}{4}G_4,
\end{eqnarray*}
which after multiplying each term by the appropriate $z(\alpha)$ becomes
\begin{eqnarray*}
&                 &   3G_1^4 + 1  G_1^2G_2 -1 G_2^2 -0 G_1G_3 - G_4  \\
&\rightarrow&   ( \chi^{(3,1)} (1^4), \chi^{(3,1)} (1^2,2), \chi^{(3,1)} (2^2), \chi^{(3,1)} (1,3), \chi^{(3,1)} (4))\\
&         =       &(3, 1, -1, 0, -1)
\end{eqnarray*}
This is the vector of values of the character at each conjugacy class as indicated. 


\vspace{0.5cm}

\noindent\textbf{P\'olya's Theory of Counting.}  The earliest use of isobaric polynomials of which we are aware is in the P\'olya theory of counting \cite{GPol,PTW}. The P\'olya theory of counting typically handles questions of the following sort: How many ways can $n$ disjoint subsets of a set be coloured with $m$ colours,  sometimes with special extra conditions.  We have a group $H$ of symmetries operating on the set and a labeling of the colours.  There are two important functions in this procedure.  One is the \emph{Cycle Indicator}: \begin{equation}
\label{ }
 C(H) = \frac{1}{|H|} \sum_{\alpha  \vdash n}  c_H(\alpha) G^{\alpha},
\end{equation}

\noindent where $G^\alpha =  G_1^{\alpha_1}G_2^{\alpha_2}\cdots G_n^{\alpha_n}$ represents the conjugacy class of the group indexed by the partition $(1^{\alpha_1}, 2^{ \alpha_2}, \ldots,n^{ \alpha_n})$ in terms of GLP,  and  where $c_H(\alpha) $ is the size of the conjugacy class indexed  by $(1^{\alpha_1}, 2^{ \alpha_2},\ldots,n^{ \alpha_n})$.  This polynomial is just the polynomial in the Burnside Theorem that counts the number of orbits of the group of symmetries.  

The other important function is the inventory function \cite[p.68]{PTW}.
The \emph{inventory} is a set of  $r$ variables $\{x_1,\ldots,x_r\},$ and the \emph{inventory function} is a function which assigns to $G_j$ the power symmetric polynomial $x_1^j+\cdots+x_r^j$. 

\bigskip
\noindent{\bf EXAMPLE.} How many ways can the vertices of a square be coloured Red and Blue if two colourings are regarded to be equal if there is a symmetry of the square that takes one into the other?  Then $n=4$, $k = 2$. The encoding of the cycle structure is given by the isobaric polynomial $$t_1^4+2t_1^2t_2+3t_2^2+2t_4.$$
The inventory is $\{x,y\}$, and the inventory function is
$$(x+y)^4+2(x+y)^2(x^2+y^2)+3(x^2+y^2)^2+2(x^4+y^4).$$
which we can write as $$G_1^4+2G_1^2G_2+3G_2^2+2G_4.$$
One can compute the total number of distinct colourings by simplifying this expression, adding the coefficients of the monomials and dividing by the order of the dihedral group of degree 4;  more simply set $x=1=y$ in the cycle index, compute and divide by the group order. The result is $48/8 = 6$;  thus there are 6 distinct colourings.  In order to compute the number of ways that a particular colouring can occur,  say, the colouring coded by $xyxy$ (defined up to a symmetry of the square), find the coefficient of the class $xyxy$ divided by $8 = |H|$.

\section{Representability}    
\label{sec:Rep} 

\begin{definition} \label{mydefinitionglobalrep}
An arithmetic function $f$ is \emph{locally representable} if for each prime $p$ there is a core polynomial, or power series $[t_1,t_2,\ldots,t_k,\ldots]$ such that $f(p^n) = F_n $ for all $n \in \mathbb{Z}$,  where $\{F_n\}$ is the GFP sequence induced by this core polynomial.
\end{definition}

\begin{definition} \label{mydefinitionlocalrep}
An arithmetic function $f$ is \emph{globally representable} if there is a core polynomial, or power series $[t_1,t_2,\ldots,t_k,\ldots]$ such that $f(n) = F_n$ for all $n\in \mathbb{Z}$,  where $\{F_n\}$ is the GFP sequence induced by this core polynomial.
\end{definition}

\begin{definition} An arithmetic function $f$ is \emph{locally linearly recursive} if for each prime $p$
$$f(p^n)=a_1f(p^{n-1})+a_2f(p^{n-2})+\cdots+a_nf(1),$$
where $a_i$'s are determined by the core polynomial for $p$. An arithmetic function $f$ is \emph{globally linearly recursive} if
$$f(n)=a_1f({n-1})+a_2f({n-2})+\cdots+a_nf(0),$$
where $a_i$'s are determined by a core polynomial. 

\end{definition}

 \begin{proposition}\label{myMFcond.}  
 A necessary and sufficient condition that $f$ be representable,  either globally or locally,  is that it be, respectively, globally or locally linearly recursive.  
\end{proposition}
 
 \proof  Since the GFP sequence is inherently a linear recursion, it is clear that the  function $f$ must also be either a locally or globally linear recursion.  On the other hand,  if it is globally or locally linearly recursive,  the parameters of the linear recursion determine a core polynomial (or power series) which in turn induces a suitable GFP sequence.  \qed
 \begin{corollary}\label{myMFLR} 
  Every multiplicative function is locally linearly recursive and hence locally representable. \qed
  \end{corollary}
  
 \begin{proposition}\cite[Rearick, Theorem 4]{R}\label{mylemmaL=0}  
Let $f $ be an arithmetic function, then $f \in \mathcal{M} $ if and only if $Lf(m)  = 0$ whenever $m$ is not a power of a prime.
 \end{proposition}

\bigskip
\noindent\textbf{REMARK.}  Let  $f  \in\mathcal{A}$ and be globally and locally representable. Since $f$ is globally representable there is a global representation $f\xrightarrow{\rm gr}\,_fF$  with $\,_f\mathcal{L}(\,_fF)=\,_fG$. Since $f$ is locally representable, $f\in\mathcal{M}$ and for any prime $p$ there is a local representation $f\xrightarrow{{\rm lr}_p}\,_fF'$ with $\,_f\mathcal{L}(\,_fF')=\,_fG'$. Then
$$\,_fG_n=\left\{\begin{array}{cl}\,_fG'_r & \hbox{ if } n=p^r, \\0 & \hbox{otherwise.}\end{array}\right.$$Since $\,_fG$ and $\,_fG'$ uniquely determine $\,_fF$ and $\,_fF'$, and $\,_fF$ and $\,_fF'$ uniquely determine $f$, $f$ is well-defined and in a certain sense only trivially globally defined. So Rearick's theorem, Proposition \ref{mylemmaL=0}, tells us that an arithmetic function $f$ can not be ``non-trivially'' both locally representable and globally representable.
 
 \vspace{0.5cm}


 \section{``Hyperbolic'' Trigonometric Functions}
\label{sec:Trig} 

\subsection{Definitions and identities}
\label{def-iden} 

Consider $F$ and $G$, and define ``hyperbolic'' SINE and ``hyperbolic'' COSINE in the formally analogous way to the standard definition for those functions in terms of $\mathcal{E}$.

\begin{definition} \label{mydefinitioncosine}
$$  C(G) =  \frac{1}{2} (\mathcal{E}(G) + \overline{\mathcal{E}(G)});$$
$$  S(G) =  \frac{1}{2} (\mathcal{E}(G) - \overline{\mathcal{E}(G)}) .$$
\end{definition}

Since $\mathcal{E}(G)=F$ and $\overline{F_n}=- {\bf t}_n$, it is easy to show that
\begin{proposition} \label{CGSG}
$$  C(G_n) =  \frac{1}{2} (F_n - {\bf t}_n);$$
$$  S(G_n) =  \frac{1}{2} (F_n +{\bf t}_n).$$
\end{proposition}\qed

Let $\delta$ be the function whose values  are $(1,0,\ldots,0,\ldots)$.

\begin{theorem}\label{mytheoremCOSSIN}
$$ C(G)^{\ast 2} - S(G)^{\ast 2} = \delta.$$
\end{theorem}

\proof \begin{eqnarray*}
&& C(G)\ast C(G)-S(G)\ast S(G) \\
& =&  \frac{1}{4}[\mathcal{E}(G)\ast \mathcal{E}(G)+\overline{\mathcal{E}(G)}\ast\overline{\mathcal{E}(G)}+2(\mathcal{E}(G)\ast\overline{\mathcal{E}(G)})]\\
&& \qquad-\frac{1}{4}[\mathcal{E}(G)\ast \mathcal{E}(G)+\overline{\mathcal{E}(G)}\ast\overline{\mathcal{E}(G)}-2(\mathcal{E}(G)\ast\overline{\mathcal{E}(G)})]\\
&=&\mathcal{E}(G)\ast\overline{\mathcal{E}(G)}\\
&=&\delta.
\end{eqnarray*}\qed

 \begin{theorem}\label{mytheoremCOS}
 Let $F$ and $G$ be induced by the core $[t_1,\ldots,t_k]$ and $F'$ and $G'$ be induced by the core $[t'_1,\ldots,t'_k]$ with $\mathcal{L}(F)=G$ and $\mathcal{L}(F')=G'$,  then
 $$C(G+G') = C(G) \ast C(G') + S(G) \ast S(G'),$$ 
 $$S(G+G' )= S(G) \ast C(G') + C(G) \ast S(G').$$
 \end{theorem}

\proof \begin{eqnarray*}
     & &  C(G+G')   \\
 & = &  \frac{1}{2} (\mathcal{E}(G+G') + \overline{\mathcal{E}(G+G')} ) \\
 & = &   \frac{1}{2}(\mathcal{E}(G)\ast \mathcal{E}(G') +\overline{\mathcal{E}(G)\ast \mathcal{E}(G') }) 
\end{eqnarray*}  
while
\begin{eqnarray*}
 &  &  C(G) \ast C(G') + S(G) \ast S(G') \\
 & = &  \frac{1}{4} [(\mathcal{E}(G)+\overline{\mathcal{E}(G)})\ast(\mathcal{E}(G')+\overline{\mathcal{E}(G')}) + (\mathcal{E}(G)-\overline{\mathcal{E}(G)})\ast (\mathcal{E}(G')-\overline{\mathcal{E}(G')})] \\
& = &  \frac{1}{2}[(\mathcal{E}(G) \ast \mathcal{E}(G') + \overline{\mathcal{E}(G)} \ast \overline{\mathcal{E}(G')}]\\
& = & \frac{1}{2}[(\mathcal{E}(G) \ast \mathcal{E}(G') + \overline{\mathcal{E}(G)\ast \mathcal{E}(G')}]. 
\end{eqnarray*}
The proof for $S$ is analogous. \qed


\subsection{Some Computations}
\label{sec:App}

We consider some applications to the ring of arithmetic functions $\mathcal{A}$. We call an arithmetic function $f$ \emph{representable} if it is either locally representable or globally representable. To illustrate Proposition \ref{CGSG}, we exhibit some computations.

Consider the multiplicative functions  $\tau$ and $\sigma$. Let $\tau\xrightarrow{{\rm lr}_p}  \, _\tau F$. Then $\, _\tau\mathcal{L}(\, _\tau F_n)=\, _\tau G_n=2$, see Section \ref{sec:LE}. Using Proposition \ref{CGSG}, it is not difficult to compute the following table for COSINE and SINE.
\begin{align*}
C(\, _\tau G_0) &=1&S(\, _\tau G_0) &= 0\\
C(\, _\tau G_1) &=0&S(\, _\tau G_1)&= 2\\
C(\, _\tau G_2) &= 2&S(\, _\tau G_2) &= 1\\
C(\, _\tau G_3) &= 2&S(\, _\tau G_3) &= 2
\end{align*}

For $\sigma$, let $\sigma\xrightarrow{{\rm lr}_p}  \, _\sigma F$. Recalling that  $\,_\sigma F_{n} = 1+p+\cdots+p^n$, we have$\,_\sigma G_{n} = p^n+1$,  for all  $n\geqslant1$, and hence the following table
\begin{align*}
C(\, _\sigma G_0) &= 1&S(\, _\sigma G_0)& = 0\\
C(\, _\sigma G_1) &=0&S(\, _\sigma G_1) &= p+1\\
C(\, _\sigma G_2) &=\frac{1}{2} (p^2+2p+1)&S(\, _\sigma G_2)& =\frac{1}{2} (p^2+1)\\
C(\, _\sigma G_3)& = \frac{1}{2} (p^3+p^2+p+1)&S(\, _\sigma G_3)& = \frac{1}{2} (p^3+p^2+p+1)\end{align*}

For the Euler totient funtion $\phi$, its core polynomial is $[p-1,\ldots]$, $\,_\phi F_{n}=p^n-p^{n-1}$ 
and $\,_\phi G_{n} = p^n-1$,  for all  $n\geqslant1$,
\begin{align*}
C(\,_\phi G_{0}) &=1&S(\,_\phi G_{0}) &= 0\\
C(\,_\phi G_{1}) &=  0&S(\,_\phi G_{1}) &= p-1\\
C(\,_\phi G_{2}) &= \frac{1}{2} (p^2-2p+1)&S(\,_\phi G_{2}) &= \frac{1}{2} (p^2-1)\\
C(\,_\phi G_{3})& = \frac{1}{2} (p^3-p^2-p+1)&S(\,_\phi G_{3}) &= \frac{1}{2} (p^3-p^2+p-1)
\end{align*}

For the zeta function $\zeta$, its core polynomial is $[1]$, $\,_\zeta F_n=1$ and 
$\,_\zeta G_{n} = 1$,  for all  $n,$
\begin{align*}
C(\,_\zeta G_{0}) &= 1&S(\,_\zeta G_{0})& = 0\\
C(\,_\zeta G_{1}) &=  0&S(\,_\zeta G_{1})& = 1\\
C(\,_\zeta G_{2})& = \frac{1}{2} &S(\,_\zeta G_{2})& = \frac{1}{2}\\
C(\,_\zeta G_{3}) &= \frac{1}{2}&S(\,_\zeta G_{3}) &= \frac{1}{2}  \\
\end{align*}

For the globally representable Catalan numbers $\Gamma(n)=\frac{1}{n+1} \binom{2n}{n}$, its core polynomial is $[1,1, 2,5,14,\ldots,t_{n+1},\ldots]$, where $t_{n+1}=\Gamma(n)$, $\,_\Gamma F_n=\Gamma(n)$ and $\,_\Gamma G_n=\,_\Gamma\mathcal{L}(\,_\Gamma F_n)={2n-1\choose n}$ for $n\geqslant 1$.
\begin{align*}
C(\,_\Gamma G_1) &=  0&S(\,_\Gamma G_1) &= 1\\
C(\,_\Gamma G_2) &= \frac{1}{2} &S(\,_\Gamma G_2) &= \frac{3}{2}\\
C(\,_\Gamma G_3) &= \frac{3}{2}&S(\,_\Gamma G_3)& = \frac{7}{2} \\
\end{align*}

\section{Acknowledgments} H. Li is supported by NSF grant DMS-0652641. 
\label{Ack}

\bigskip
\noindent Huilan Li

\noindent Department of Mathematics, Drexel University, Philadelphia PA 19104 U.S.A

\noindent Email: \emph{huilan.li@gmail.com}

\smallskip
\noindent Trueman MacHenry 

\noindent Department of Mathematics and Statistics, Toronto, Ontario M3J 1P3 CANADA

\noindent Email: \emph{machenry@mathstat.yorku.ca}

\smallskip
\noindent\today


\begin{thebibliography}{\hphantom{99}}

\bibitem{CG} T. B. Carroll and A. A. Gioia,  {On a subgroup of the group of multiplicative arithmetic functions}, J. Austral. Math. Soc. Ser. A \textbf{20} (1975), no.3, 348--358.

\bibitem{CE} E. D. Cashwell and C. J. Everett,  {The ring of number theoretic functions}, Pacific J. Math. \textbf{9 }(1959), no.4, 975--985.

\bibitem{JE} J. Elliott, {Ring structures on groups of arithmetic functions}, J. of Number Theory, \textbf{128}, (2008) 709--730.

\bibitem{online} The On-Line Encyclopedia of Integer Sequences.

\bibitem{L} D. H. Lehmer, {An extended theory of Lucas' functions}, Annals of Mathematics, 2nd Ser., \textbf{31} (1930), 419--448.

\bibitem{Macd}  I. G. Macdonald,  {Symmetric Functions and Hall Polynomials}, Clarendon Press, Oxford, (1995).
 
\bibitem{TM}  T. MacHenry, {A Subgroup of the Group of Units in the Ring of Arithmetic Functions}, Rocky Mountain Journal of Mathematics, \textbf{29}, (1999), 1055--1065.

\bibitem{TM2}  T. MacHenry, {Generalized Fibonacci and Lucas Polynomials and Multiplicative Arithmetic Functions}, Fibonacci Quarterly, \textbf{38}, (2000), 17--24.

\bibitem{MT}  T. MacHenry  and G. Tudose, {Reflections on Isobaric Polynomials and Arithmetic Functions}, Rocky Mountain Journal of Mathematics, Vol. \textbf{35}, No.3, (2005), 901--928.

 \bibitem{MT2}  T. MacHenry and G. Tudose, {Differential Operators and Weighted Isobaric Polynomials}, Rocky Mountain Journal of Mathematics, Vol \textbf{36} N6 (2006) 1957--1976.

\bibitem{MW} T.  MacHenry and K. Wong, {A Correspondence Between Isobaric Rings and Multiplicative Arithmetic Funtions},  Rocky Mountain Journal of Mathematics, to appear. 

\bibitem{MW2} T. MacHenry and K. Wong, {Degree-$k$ linear recursions modulo ($p$) and algebraic number fields}, Rocky Mountain Journal of Mathematics, to appear (2010).

\bibitem{PJM}P. J. McCarthy, {Introduction to Arithmetical Functions}, Unversitext, Springer, New York, (1986).

\bibitem{GPol} G. P\'olya, {Kombiinatorische Anzahlbestimmungen f\"ur Gr\"uppen}, Graphen und chemische Verbindungen, Acta  Math., \textbf{68}, (1937), 145--254.

\bibitem{PTW} G. P\'olya, R. E. Tarjan and D. R. Woods, {Notes on Introductory Combinatorics}, Birkh\"auser, Progress in Computer Science, (1983).

\bibitem{R} D. Rearick, {Operators on algebras of arithmetic functions}, Duke Math. J. \textbf{35} (1968), no. 4, 761--766.

\bibitem{R2} D. Rearick, {The trigonometry of numbers}, Duke Math. J. \textbf{35} (1968), no. 4, 767--776.

\bibitem{ReS} D. Redman and R. Sivaramakrishnan, {Some properties of specially multiplicative functions}, J. Number Theory, \textbf{13} (1981), no.2, 210--227.

\bibitem{JR} J. Rutkowski, {On recurrence characterization of rational arithmetic functions}, Funct. Approx. Comment. Math. \textbf{9} (1980), 45--47.

\bibitem{HS} H. N. Shapiro,  {On the Convolution Ring of Arithmetic Functions}, Communications of Pure and Applied Mathematics, vol. xxv, (1972), 287--336. 


\end{thebibliography}
\end{document}